\documentclass[11pt,twoside]{gsp}

\def\semicolon{\nobreak\mskip2mu\mathpunct{}\nonscript\mkern-\thinmuskip{;}
\mskip6muplus1mu\relax} 

\def\Id{
{\rm Id}\, }

\newcommand{\e}{\mathrm{e}}

\renewcommand{\d}{\mathrm{d}}

\theoremstyle{plain}
\newtheorem{theorem}{Theorem}
\newtheorem{definition}[theorem]{Definition}

\newtheorem{example}[theorem]{Example}

\newtheorem{remark}[theorem]{Remark}

\newenvironment{proof}[1][{}]{\par\noindent\textbf{Proof{#1}: }}{\hspace*{\fill}$\blacksquare$\smallskip\noindent\par}

\begin{document}

\thispagestyle{plain}

\title[Commuting Pairs of Generalized Contact Metric Structures]{Commuting Pairs of Generalized Contact Metric Structures}
\author[Janet Talvacchia]{Janet Talvacchia}

\date{February, 27, 2018}

\maketitle

\comm{Communicated by XXX}\

\begin{abstract}
In this paper, we prove a theorem that gives a simple criterion for generating commuting pairs of generalized almost complex structures on spaces that are the product of two generalized almost contact metric spaces.  We examine the implications of this theorem with regard to the definitions of generalized Sasakian and generalized co-K\"ahler geometry.
\\[0.2cm]
\end{abstract}

\label{first}

\section{Introduction}\label{S:intro}
\indent
\indent The notion of a generalized complex structure, introduced by Hitchin in his paper \cite{H} and developed by Gualtieri \cite{G1, G2} is a framework that unifies both complex and symplectic structures.  These structures exist only on even dimensional manifolds. The odd dimensional analog of this structure, a generalized contact structure, was taken up by Vaisman \cite {V1, V2},  Poon, Wade \cite{PW}, Sekiya \cite{S}, and Aldi and Grandini \cite{AG}. This framework unifies almost contact, contact, and cosymplectic structures.  Generalized K\"ahler structures were introduced by Gualtieri \cite{G1,G2, G3} and have already found their way into the physics literature \cite{Hu, LMTZ, Gr}.  A natural question to ask is when the product of two generalized almost contact metric structures produce a generalized K\"ahler structure.

\indent  K\"ahler structures are commuting pairs of generalized complex structures whose product forms a generalized metric.  In this article, we prove the following  theorem that gives a simple criterion for generating commuting pairs of generalized almost complex structures on spaces that are the product of generalized almost contact metric spaces.  This reduces the assessment of whether the spaces could be generalized K\"aher to integrability issues.
\begin{theorem}\label{T1}
Let $M_{1}$ and $M_{2}$ be odd dimensional smooth manifolds.
Assume that $M_1$  and $M_2$ each have two generalized almost contact metric structures. Denote these structures by
$( \Phi_1, E_{\pm_1}, G_1)$  and  $(\tilde{\Phi}_1,  \tilde{E}_{\pm_1}, G_1)$ on $M_1$, and $( \Phi_2, E_{\pm_2}, G_2)$  and  $(\tilde{\Phi}_2 , \tilde{E}_{\pm_2}, G_2)$ on $M_2$.  Let $\mathcal{J}_1$  be  the natural generalized complex structure on $M_1 \times M_2$ induced from $\Phi_1$ and $\Phi_2$  and let  $\mathcal{J}_2$ be the natural generalized almost complex structure on $M_1 \times M_2$ induced from $\tilde{\Phi}_1$ and $\tilde{\Phi}_2$.  Then $\mathcal{J}_1$ and $\mathcal{J}_2$ form a commuting pair of  generalized almost complex structures if and only  if $\Phi_i$ commutes with $\tilde{\Phi}_i$ for $ i=1,2$  and if either  the condition $E_{\pm 1 }= \tilde{E}_{\pm 1}$ and $E_{\pm 2}= \tilde{E}_{\pm 2}$ or the condition $E_{\pm 1} = \tilde{E}_{\mp 1}$ and $E_{\pm 2} = \tilde{E}_{\mp 2}$ holds.
\end{theorem}

\indent  The value of this theorem is twofold.  It gives a simple and computable criterion for determining, when handed two generalized almost contact metric spaces, if their product is possibly generalized K\"ahler directly from data on the generalized almost contact metric spaces. One can compute the relevant commuting pairs of commuting generalized almost complex structures and examine their integrability. This facilitates analyzing generalized K\"aher structures arising from classical Sasakian manifolds. (See Examples \ref{ex17} and \ref{ex18} below.) Secondly it suggests that notions of generalized co-K\"ahler and generalized Sasakian can be reformulated in terms of commuting pairs of generalized almost contact structures.  This adds clarity and elegance to this topic since generalized K\"ahler structures are defined in terms of commuting pairs of generalized complex structures and it highlights the subtle role of the integrability criteria.

\indent In Section 2 we gather the basics of  generalized complex, generalized K\"ahler, and generalized contact geometry. In Section 3, we prove the main theorem stated above. Then, in Section 4, we examine the implications of this theorem for the definitions of generalized co-K\"ahler and generalized Sasakian geometry.

\section{Preliminaries}\label{B:Back}

\indent We begin with a very brief review of generalized geometric structures. Throughout this paper we let $M$ be a smooth manifold. Consider the big tangent bundle, $TM\oplus ~T^*M$.  We define a neutral metric on $TM\oplus T^*M$ by$$  \langle X + \alpha  , Y + \beta \rangle =  \frac{1}{2} (\beta (X) + \alpha (Y) )$$ and the Courant bracket by $$[[X+\alpha, Y+ \beta ]] = [X,Y] + {\mathcal L}_X\beta -{\mathcal L}_Y\alpha -\frac{1}{2} d(\iota_X\beta - \iota_Y\alpha)$$ where $X, Y \in TM$ and $\alpha ,\beta  \in T^*M$ . A subbundle of $TM\oplus T^*M$ is said to be involutive  or integrable if its sections are closed under the Courant bracket\cite{G1}.

\begin{definition}
A generalized  almost complex structure on $M$ is an endomorphism $\mathcal J$ of $TM\oplus T^*M$ such that $\mathcal J + \mathcal J^* =  0 $ and $\mathcal J^2 = - \Id$. If the $\sqrt{-1}$ eigenbundle $L\subset (TM\oplus TM^{*})\otimes\mathbb{C}$ associated to $\mathcal J$ is involutive with respect to the
Courant bracket, then $\mathcal J$ is called a generalized complex structure.
\end{definition}

\indent Here are the prototypical examples:

\begin{example}[\!\cite {G1}]\hspace*{-0.12cm}{\bf.}\label{ex3}
Let $(M^{2n}, J)$ be a complex structure.  Then we get a generalized complex structure by setting
$$\mathcal J_{J} = \left ( \begin{array}{rc}  -J & 0 \\ 0 & J^* \end{array} \right ).$$
\end{example}

\begin{example}[\!\cite {G1}]\hspace*{-0.12cm}{\bf.}\label{ex4}
Let $(M^{2n}, \omega )$ be a symplectic structure.  Then we get a generalized  complex structure by setting
$$\mathcal J_{\omega} = \left ( \begin{array}{cc}  0 & -\omega^{-1} \\ \omega & 0 \end{array} \right ).$$
\end{example}
\indent Diffeomorphisms of $M$ preserve the Lie bracket of smooth vector fields and in fact such diffeomorphisms are the
only automorphisms of the tangent bundle. But in generalized geometry, there is actually more flexibility.
Given $T\oplus T^{*}$ equipped with the Courant bracket, the automorphism group is comprised of the diffeomorphisms of $M$
and some additional symmetries called \emph{B}-field transformations \cite{G1}.
\begin{definition}[\!\cite {G1}]\hspace*{-0.12cm}{\bf.}
Let $B$ be a closed two-form which we view as a map from $T \rightarrow T^{*}$ given by interior product. Then the invertible bundle map
$$\e^{B}:= \left ( \begin{array}{cc}  1 & 0 \\ B & 1 \end{array} \right):X+\xi \longmapsto X+\xi + \iota_{X}B$$
is called a B-field transformation.
\end{definition}
A $B$-field transformation of a generalized (almost) complex structure $( M, \e^B {\mathcal J} \e^{-B})$ is again a  generalized (almost) complex structure.

\indent Recall that we can reduce the structure group of $T\oplus T^{*}$ from $O(2n,2n)$ to the maximal compact subgroup $O(2n)\times O(2n)$. This
is equivalent to an orthogonal splitting of $T\oplus T^{*}=V_{+}\oplus V_{-}$, where $V_{+}$ and
$V_{-}$ are positive and negative definite respectfully with respect to the inner product. Thus we can define a positive definite Riemannian metric
on the big tangent bundle by $$G=\langle\,,\rangle|_{V_{+}}-\langle\,,\rangle|_{V_{-}}.$$
A positive definite metric $G$ on $M$ is an automorphism of $TM\oplus T^*M$ such that $G^{*}=G$ and $G^{2}=1.$ In the presence
of a generalized almost complex structure $\mathcal J_1$, if $G$ commutes with $\mathcal J_1$ ( $G\mathcal J_1 = \mathcal J_1  G$) then $G\mathcal J_1 $ squares to $-1$ and we generate a second generalized almost complex structure, $\mathcal J_2$ $= G\mathcal J_1$, such that $\mathcal J_1$ and $\mathcal J_2$ commute and $G=-\mathcal J_1 \mathcal J_2$.
This leads us to the following
\begin{definition}[\!\cite {G1}]\hspace*{-0.12cm}{\bf.} A  generalized K\"ahler structure is a pair of commuting generalized complex structures $\mathcal J_{1}, \mathcal J_{2}$ such that
$G=-\mathcal J_{1}\mathcal J_{2}$ is a positive definite metric on $T\oplus T^{*}.$
\end{definition}

\begin{example}[\!\cite {G1}]\hspace*{-0.12cm}{\bf.}
 Consider a K\"ahler structure $(\omega,J,g)$ on $M$. By defining $\mathcal J_{J}$ and $\mathcal J_{\omega}$ as in
{\rm Examples \ref{ex3}} and {\rm\ref{ex4}}, we obtain a generalized K\"ahler structure on $M$, where
$$G=\left ( \begin{array}{cc}  0 & g^{-1} \\ g & 0 \end{array} \right ).$$

\end{example}
\indent We now recall the odd dimensional analog of generalized complex geometry. We use the definition given by Sekiya  (see \cite{S}).
\begin{definition} A generalized almost contact structure on $M$ is a triple $(\Phi, E_\pm)$ where $\Phi $ is an endomorphism of $TM\oplus T^*M$, and $E_+$ and $E_-$ are sections of $TM\oplus T^*M$ which satisfy
\begin{equation}\label{phi}
\Phi + \Phi^{*}=0, \qquad
\Phi \circ \Phi = -\Id + E_+ \otimes E_- + E_- \otimes E_+
\end{equation}
\begin{equation}\label{sections}
 \langle E_\pm, E_\pm \rangle = 0,  \qquad     2\langle E_+, E_-  \rangle = 1.
\end{equation}

\end{definition}
An easy and immediate consequence of these definitions (see \cite{S}) is
\begin{equation}\label{PhivanishE}
\Phi(E_{\pm})=0.
\end{equation}
Now, since $\Phi$ satisfies $\Phi^3 + \Phi =0$, we see that $\Phi$ has $0$ as well as $\pm \sqrt{-1}$ eigenvalues when viewed as an endomorphism of the complexified big tangent bundle $(TM\oplus ~T^*M) ~\otimes { ~\mathbb C}$.  The kernel of $\Phi$ is $L_{E_+} \oplus L_{E_-}$ where $L_{E_\pm}$ is the line bundle spanned by ${E_\pm}$.  Let $E^{(1,0)}$ be the $\sqrt{-1}$ eigenbundle of $\Phi$.  Let $E^{(0,1)}$ be the $-\sqrt{-1}$ eigenbundle. Observe

$$
E^{(1,0)} = \lbrace X + \alpha - \sqrt{-1}  \Phi (X + \alpha ) \semicolon  \langle E_\pm, X + \alpha \rangle = 0 \rbrace
$$

$$
E^{(0,1)} = \lbrace X + \alpha + \sqrt{-1}  \Phi (X + \alpha ) \semicolon \langle E_\pm, X + \alpha \rangle = 0 \rbrace .$$

Then the complex vector bundles
$$L^+ = L_{E_+} \oplus E^{(1,0)}$$
and
$$L^- = L_{E_-} \oplus E^{(1,0)}$$
are maximal isotropics.
\begin{definition}
A generalized almost contact structure $(\Phi,E_{\pm})$ is a  generalized contact structure if either $L^{+}$ or $L^{-}$ is closed with respect to the Courant bracket. The generalized
contact structure is strong if both $L^{+}$ and $L^{-}$ are closed with respect to the Courant bracket.
\end{definition}

\begin{definition}
A  generalized almost contact structure $(M, \Phi,E_{\pm})$ is a normal generalized contact structure if  $\Phi$ is strong and $[[E_+, E_-] ] = 0$.
\end{definition}

\begin{remark}  This definition of normality is motivated by {\rm Theorem 1} of {\rm\cite{GT1}} that shows that product of two generalized almost contact spaces $(M_1, \Phi_1, E_{\pm 1})$  and $(M_2, \Phi_2, E_{\pm 2})$ induces a natural generalized almost  complex structure on $M_1 \times M_2 $.  The generalized complex structure is integrable if each $\Phi_i$ is strong and  $[[E_{+i},  E_{-i}]] = 0$.
\end{remark}
Here are the standard examples

\begin{example}[\!\cite {PW}]\hspace*{-0.12cm}{\bf.}\label{almost contact example}
Let $(\phi , \xi, \eta)$ be a normal almost contact structure on a manifold $M^{2n+1}$.  Then we get a generalized almost contact structure by setting
$$ \Phi = \left ( \begin{array}{cc}  \phi & 0 \\ 0 & -\phi^* \end{array} \right ),\qquad  E_+ = \xi,\qquad  E_-= \eta $$  where $(\phi^*\alpha )(X) = \alpha (\phi (X)), \   X \in TM,\   \alpha \in T^{*}M$. Moreover, $(\Phi, E_\pm)$ is an example of a strong generalized almost contact structure.
\end{example}

\begin{example}[\!\cite {PW}]\hspace*{-0.12cm}{\bf.} \label{contact example}
Let $( M^{2n+1}, \eta )$ be a contact manifold with $\xi $ the corresponding Reeb vector field so that
$$ \iota_\xi \d \eta = 0, \qquad  \eta ( \xi ) = 1.$$
Then $$\rho ( X) := \iota_X \d \eta - \eta ( X)\eta$$ is an isomorphism from the tangent bundle to the cotangent bundle.  Define a bivector field by
$$\pi (\alpha , \beta ) := \d \eta (\rho^{-1}( \alpha ), \rho^{-1}( \beta ))$$
where $\alpha, \beta \in T^{*}. $
We obtain a generalized almost contact structure by setting
$$ \Phi = \left ( \begin{array}{cc}  0 & \pi \\ \d \eta & 0 \end{array} \right ),\qquad   E_+ = \eta,\qquad   E_-= \xi .$$
In fact, $(\Phi, E_{\pm})$ is an example which is not strong.
\end{example}

Finally we have
\begin{definition}

A generalized almost contact metric structure is a generalized almost contact structure
$(\Phi , E_{\pm})$ along with a generalized Riemannian metric $G$ that satisfies
\begin{equation}\label{compatcond}
-\Phi G \Phi = G - E_+ \otimes E_+ -E_- \otimes E_-.
\end{equation}

\end{definition}

\section{ Proof of the Main Theorem}

In this section we give the proof of Theorem 1 which was stated in the introduction.

\begin{proof}
Since $( M_i, \Phi_i, E_{\pm_i}, G_i)$ are generalized almost contact metric structures for $i=\{1,2\}$, $M_1 \times M_2$ with product metric $G=G_1+G_2$ is a generalized almost complex manifold (see \cite{GT1}) and we can form the generalized almost complex structures $\mathcal{J}_1$  and $\mathcal{J}_2$ as follows
\begin{eqnarray}\label{3.1}
\mathcal{ J}_1 (X_1 +\alpha_1,X_2 +\alpha_2 ) &= &( \Phi_1(X_1+\alpha_1)- 2\langle E_{+2},X_2 +\alpha_2 \rangle E_{+1}\nonumber\\
 && - 2\langle E_{-2}, X_2 +\alpha_2 \rangle E_{-1},\nonumber\\[-0.2cm]
 \\[-0.3cm]
 &&\Phi_2(X_2 +\alpha_2)  + 2\langle E_{+1}, X_1 +\alpha_1 \rangle E_{+2}\nonumber \\
 &&+ 2\langle E_{-1}, X_1 +\alpha_1 \rangle E_{-2} ).\nonumber
\end{eqnarray}

\begin{eqnarray}\label{3.2}
\mathcal{J}_2 (X_1 +\alpha_1,X_2 +\alpha_2 ) &=& ( \tilde{\Phi}_1(X_1+\alpha_1)- 2\langle \tilde{E}_{+2},X_2 +\alpha_2 \rangle \tilde{E}_{+1}\nonumber\\
 && - 2\langle \tilde{E}_{-2}, X_2 +\alpha_2 \rangle \tilde{E}_{-1}, \nonumber\\[-0.2cm]
 \\[-0.3cm]
&&\tilde{\Phi}_2(X_2 +\alpha_2)  + 2\langle \tilde{E}_{+1}, X_1 +\alpha_1 \rangle \tilde{E}_{+2}\nonumber\\
 &&+ 2\langle \tilde{E}_{-1}, X_1 +\alpha_1 \rangle \tilde{E}_{-2} ). \nonumber
\end{eqnarray}

By direct computation we see
\begin{eqnarray}\label{3.3}
\mathcal{J}_1\mathcal{J}_2 &=& (\Phi_1 \tilde{\Phi}_1  - 2\langle E_{-1}, X_1 + \alpha_1\rangle E_{+1} -2 \langle E_{+1}, X_1 + \alpha_1 \rangle E_{-1}, \nonumber\\[-0.2cm]
\\[-0.3cm]
&&\Phi_2 \tilde{\Phi}_2  -2\langle E_{-2}, X_2 + \alpha_2 \rangle E_{+2} -2 \langle E_{+2}, X_2 +\alpha_2\rangle E_{-2})\nonumber
\end{eqnarray}

\begin{eqnarray}\label{3.4}
\mathcal{J}_2\mathcal{J}_1 &=& (\tilde{\Phi}_1 \Phi_1  - 2\langle \tilde{E}_{-1}, X_1 + \alpha_1\rangle \tilde{E}_{+1} -2 \langle \tilde{E}_{+1}, X_1 + \alpha_1 \rangle \tilde{E}_{-1}, \nonumber\\[-0.2cm]
\\[-0.3cm]
&&\tilde{\Phi}_2 \Phi_2 -2\langle \tilde{E}_{-2}, X_2 + \alpha_2 \rangle \tilde{E}_{+2} -2 \langle \tilde{E}_{+2}, X_2 +\alpha_2\rangle \tilde{E}_{-2}).\nonumber
\end{eqnarray}

So we see $\mathcal{J}_1\mathcal{J}_2 = \mathcal{J}_2\mathcal{J}_1 $ if and only if $\Phi_i \tilde{\Phi}_i = \tilde{\Phi}_i \Phi_i$ for $i=\{1,2\}$ and either  the condition $E_{\pm 1 }= \tilde{E}_{\pm 1}$ and $E_{\pm 2}= \tilde{E}_{\pm 2}$ or the condition $E_{\pm 1} = \tilde{E}_{\mp 1}$ and $E_{\pm 2} = \tilde{E}_{\mp 2}$ holds.

\end{proof}

\begin{remark}
The integrability of $\mathcal{J}_1$ and $\mathcal{J}_2$ can be determined via {\rm Theorem 1} in {\rm\cite{GT1}}.
\end{remark}

$M_1 \times M_2$ may also have commuting generalized almost complex structures $\tilde{\mathcal{J}}_1$ and $\tilde{\mathcal{J}}_2$ that have nothing to do with the underlying generalized almost contact structures $\Phi_1, \tilde{\Phi}_1, \Phi_2, \tilde{\Phi}_2$ on $M_1 $  and $M_2$ and, thus, $\tilde{\mathcal{J}}_1$ and $\tilde{\mathcal{J}}_2$ would not have the same algebraic form as $\mathcal{J}_1$ and $\mathcal{J}_2$. So we emphasize that the ``only if'' part of the above theorem is restricted to generalized almost complex structures on $M_1 \times M_2$ that are induced from the underlying generalized almost contact structures on $M_1$ and $M_2$.

\section{Implications for Generalized Sasakian and Generalized Co-K\"ahler Manifolds}

In a previous paper (\!\!\cite{GT2}, joint work with Ralph Gomez),  we considered what we termed generalized co-K\"ahler structures since each classical co-K\"ahler manifold falls into this category.   A generalized co-K\"ahler structure,  $( M , \Phi , E_{\pm}, G)$, is a normal generalized contact metric structure such that both $\Phi$ and $G\Phi$ are strong. The result we proved there is as follows
\begin{theorem}\label{T4.1}

Let $M_{1}$ and $M_{2}$ be odd dimensional smooth manifolds each with a generalized contact metric structure
$(\Phi_i,E_{\pm,i},G_i),i=1,2$ such that on the product $M_1\times M_2$ are two
generalized almost complex structures: $\mathcal{J}_1$ which is the natural generalized almost complex structure induced from $\Phi_1$ and $\Phi_2$ and $\mathcal{J}_2=G\mathcal{J}_1$ where
$G=G_1\times G_2$. Then $(M_1\times M_2,\mathcal{J}_1,\mathcal{J}_2)$ is generalized K\"ahler if and only if $(\Phi_{i},E_{\pm,i},G_i)$, $i=1,2$ are  generalized co-K\"ahler structures.

\end{theorem}
Here, the commuting pair on each $M_i$ is $\Phi_i $ and $G_i\Phi_i$.  The natural generalized complex structures $\mathcal{J}_1$ formed from $\Phi_1$ and  $\Phi_2 $  is  Courant integrable  since  the $\Phi_i$  are normal. $\mathcal{J}_2 = G\mathcal{J}_1$ is the same as the generalized almost complex structure formed from $G_1\Phi_1$ and  $G_2\Phi_2 $.  $\mathcal{J}_2$ is Courant integrable since $G_1\Phi_1$ and  $G_2\Phi_2 $ are normal generalized contact metric structures.  
\par

Classical Sasakian manifolds are always outside of the generalized co-K\"ahler case. A Sasakian manifold $M$ is a normal contact metric structure $( M, \phi , \eta , \xi, g)$ where $\phi$ is a $(1,1)$ tensor field, $\xi$ is a vector field and $\eta$ is a one-form, given by the following conditions
\begin{center}
$\phi^{2}=-I+\eta\otimes\eta, \hspace{.5cm}\eta(\xi)=1$
\end{center}
and where $g$ is a Riemannian metric subject to the following compatibility condition
\begin{center}
$g(\phi X, \phi Y)=g(X,Y)-\eta(X)\eta(Y)$
\end{center}
for an vector fields $X,Y\in \Gamma(TM)$.  If we form the generalized contact structures related to the Sasakian structure we get a  strong generalized contact structure associated to the classical almost contact metric structure
$$ \Phi_\phi = \left ( \begin{array}{cc}  \phi & 0 \\ 0 & -\phi^* \end{array} \right ),\qquad  E_+ = \xi,\qquad   E_-= \eta. $$  The compatible generalized metric is $$G=\left ( \begin{array}{cc}  0 & g^{-1} \\ g & 0 \end{array} \right ).$$

However
$$G\Phi = \Phi_\eta = \left ( \begin{array}{cc}  0 & \pi \\ \d \eta & 0 \end{array} \right ),\qquad   E_+ = \eta,\qquad    E_-= \xi $$
is never strong. (See Examples \ref{almost contact example} and \ref{contact example}.)  So classical Sasakian spaces are not generalized co-K\"ahler but Theorem \ref{T1} implies that the generalized almost contact structures arising from the classical almost contact structure and the classical contact structure generate commuting pairs of generalized almost complex structures.  The issue is integrability.
\begin{example}\label{ex17}
Let $ M_1 = ( M, \phi , \eta , \xi, g)$ be a classical Sasakian space. Let
\begin{center}
$\Phi_1 =  \Phi_\phi = \left ( \begin{array}{cc}  \phi & 0 \\ 0 & -\phi^* \end{array} \right ),\qquad   E_{+ 1}= \xi,\qquad   E_{-1}= \eta $
\end{center}
  be the generalized contact structure generated  by the classical almost contact structure $\phi$. The compatible generalized metric is $$G_1 =\left ( \begin{array}{cc}  0 & g^{-1} \\ g & 0 \end{array} \right ).$$
Let

$$\tilde{\Phi}_1 = G\Phi = \Phi_\eta = \left ( \begin{array}{cc}  0 & \pi \\ \d \eta & 0 \end{array} \right ),\qquad   E_{+1} = \eta,\qquad   E_{-1}= \xi ,   \qquad  G _1$$
 be the second commuting generalized contact metric structure on $M_1$.  Let
 $$M_2 = \mathbb{R}^{+}, \qquad   \Phi_2 = 0, \qquad  E_{+2} = (0, \d t),\qquad  E_{-2} = (\frac{\partial}{\partial t}, 0)$$
 be a generalized contact metric structure on $\mathbb{R}^+$ with compatible metric
\begin{center}
$G_2 =  \left ( \begin{array}{cc}  0 & (\d t^2)^{-1} \\ \d t^2 & 0 \end{array} \right )$.
\end{center}
Let

$$\tilde{\Phi}_2  = G_2\Phi_2 = 0, \qquad \tilde{E}_{+2} = ( \frac{\partial}{\partial t} , 0),\qquad  \tilde{E}_{+2} =(0, \d t), \qquad  G_2$$

be the second commuting generalized contact metric structure on $\mathbb{R}^+\!.$
Let $\mathcal{J}_1$ be the generalized  complex structure formed from $\Phi_1$ and $\Phi_2$.  It is Courant integrable since each $\Phi_i$ is normal.  $\mathcal{J}_2$, the generalized almost complex structure formed from $\tilde{\Phi}_1$ and $\tilde{\Phi}_2$  satisfies $\mathcal{J}_2 = G\mathcal{J}_1$ but  $\mathcal{J}_2$  is not Courant integrable since $\tilde{\Phi}_1$ is not.
\end{example}
On the other hand we know, given a classical Sasakian manifold $M$, there is a generalized K\"ahler structure induced on $M\times \mathbb{R}^+$ from the classical K\"ahler structure $(M \times \mathbb{R}^+, \omega = \d (\e^t\eta), J_\phi, \e^{2t}(g+\d t^2))$ on the cone.  We construct in detail the commuting pairs of generalized almost contact metric structures on $M$ that yield the generalized complex structures arising from the classical K\"ahler structure

\begin{example}\label{ex18}
Given a classical Sasakian manifold  $( M, \phi , \eta , \xi, g)$, form the following commuting pair of generalized contact structures on $M\equiv M_1$. Let $R$ be the endomorphism of $TM_1 \oplus T^*M_1$

$$R =  \left ( \begin{array}{cc}  \e^{-t} & 0 \\ 0 & \e^t \end{array} \right ) \qquad for \qquad t\in\mathbb{R}^+.$$
 Define
 $$\Phi_1 = R\Phi_\phi R^{-1} = \Phi_\phi,\qquad
E_{+1}= R( \xi, 0),\qquad  E_{-1} = R(0, \eta)$$
   and
  $$G_1 = R \left ( \begin{array}{cc}  0 & g^{-1} \\ g & 0 \end{array} \right ) R^{-1}.$$

   Then $( M_1, \Phi_1, E_{\pm 1}, G_1)$ is a generalized contact metric structure. 
   
  Define
  
 $$ \tilde{\Phi}_1 = R\Phi_\eta R^{-1} =   RG\Phi_\phi R^{-1} = RGR^{-1}\Phi_\phi$$
 $$ \tilde{E}_{+1} = R(0, \eta), \qquad  \tilde{E}_{-1}= R( \xi, 0).$$

  Then $( M_1,\tilde{ \Phi}_1, \tilde{E}_{\pm 1}, G_1)$ is a second generalized contact metric structure, $\Phi_1$ commutes with $\tilde{\Phi}_1$ and  $E_{\pm 1} = \tilde{E}_{\mp 1}$ holds.

 \noindent  For the second factor,  let
  $$M_2 = \mathbb{R}^{+}, \qquad  \Phi_2 = 0, \qquad E_{+2} = R(0, \d t)\qquad
  E_{-2} =  R(\frac{\partial}{\partial t}, 0)$$ and $$G_2 = R \left ( \begin{array}{cc}  0 & (\d t^2)^{-1} \\ \d t^2 & 0 \end{array} \right )R^{-1}.$$ Then $( M_2, \Phi_2, E_{\pm 2}, G_2)$ is a generalized contact metric structure.

 A second  generalized contact metric structure on $\mathbb{R}^+$ satisfying the necessary conditions is obtained by  defining
  $$\tilde{\Phi}_2  = RG_2\Phi_2R^{-1} =  RG_2R^{-1}\Phi_2 = 0$$
  $$ \qquad \tilde{E}_{+2} =~ R(0, \d t),  \tilde{E}_{-2} =R ( \frac{\partial}{\partial t} , 0)$$

 Computing  according to formulas {\rm(\ref{3.1})} and {\rm(\ref{3.2})}, we see

\begin{align}
\mathcal{J}_1 ( X_1 + \alpha_1, X_2 + \alpha_2)  = &(\Phi_\phi -2\langle \d t , X_2 + \alpha _2 \rangle  \xi  - 2\langle  \frac{\partial}{\partial t}, X_2 + \alpha_2  \rangle  \eta ,\nonumber  \\[-0.3cm]
\\[-0.3cm]
 &  2\langle  \xi, X_1+ \alpha_1 \rangle \d t + 2\langle  \eta , X_1 + \alpha_1\rangle \frac{\partial}{\partial t})\nonumber
\end{align}

\begin{align}
\mathcal{J}_2 ( X_1 + \alpha_1, X_2 + \alpha_2 ) = &(R\Phi_\eta R^{-1} - 2e^{2t} \langle {\frac{\partial}{\partial t}}, X_2 + \alpha _2 \rangle    \eta - 2e^{-2t} \langle  \d t , X_2 + \alpha_2 \rangle  \xi , \nonumber\\[-0.2cm]
\\[-0.3cm]
&  2e^{2t}\langle \eta , X_1 + \alpha_1 \rangle  \d t + 2e^{-2t}\langle   \xi , X_1 + \alpha _1 \rangle  \frac{\partial }{\partial t}).\nonumber
\end{align}
\end{example}

Note that if one starts with the classical normal almost contact structure $\phi$ on $M$, computes the corresponding classical complex structure $J$ on the cone, and then computes the generalized complex structure corresponding to this classical complex structure what one gets is precisely $\mathcal{J}_1$. Similarly, if one starts with the classical contact structure, lifts to the standard symplectic structure $\omega = \d (\e^t\eta)$ on the cone, and then constructs the generalized complex structure corresponding to this classical structure, what one gets is precisely $\mathcal{J}_2$.   Hence we recover the generalized K\"ahler structure  $(M_1 \times ~\mathbb{R}^+, \mathcal{J}_1, \mathcal{J}_2)$ corresponding to the classical K\"ahler structure on $M\times \mathbb{R}^+$ from this family of commuting pairs of generalized almost contact structures. It is interesting to note that $\mathcal{J}_2 = RGR^{-1}\mathcal{J}_1 = RG\mathcal{J}_1R^{-1}$ where $G$ is the standard product metric constructed from the Riemanian metrics on $M$ and $\mathbb{R}^+$.  That is, $\mathcal{J}_2$ can be obtained as well by applying the warped metric $RGR^{-1}$ to $\mathcal{J}_1$.  The fact that $\mathcal{J}_2$ is also equal to $RG\mathcal{J}_1R^{-1}$ will be useful when we compare defining generalized Sasakian via commuting pairs of generalized almost contact to previous definitions of  generalized Sasakian below.

Examples \ref{ex17} and \ref{ex18} also point out an interesting phenomena.   In Example \ref{ex17}, $\mathcal{J}_2$ is not integrable with respect to Courant bracket but is with respect to the Wade bracket (see \cite{IW}).  Specifically Ingleas-Ponte and Wade show in (\!\cite{IW}) that, in general, generalized almost contact structures arising from either classical almost contact structures or classical contact structures are closed  as $\mathcal{E}^1(M)$ Dirac structures under the Wade bracket.  They also show a one to one correspondence between generalized almost contact structure integrable with respect to the Wade bracket and operators $\mathcal{J}$ on $M\times \mathbb{R}^+$, closed under the Wade bracket, such that $\mathcal{J}^2 = -1$.  In Example \ref{ex18},  the example corresponding to the classical K\'ahler structure on $M\times \mathbb{R}^+$, what we have is really a family of  generalized almost contact structures on $M$ parameterized by $t$.   For each $t$, $ \tilde{\Phi}_1$ is not closed with respect to Courant bracket, but it is closed with respect to the derived bracket generated by $D= \e^t\d \e^{-t} $  (as are $\Phi_1, \Phi_2$ and $\tilde{\Phi}_2$, see \cite{K, R}).  The  commuting almost complex structures generated by these objects are Courant integrable.  The fact that it is a family of generalized almost contact structures combining to form the generalized complex structure as opposed to a single generalized almost contact metric structure on $M$  producted with one on $\mathbb{R}^+$ makes this example not in contradiction to Theorem \ref{T4.1}.  These examples suggest that generalized K\"ahler structures on products of generalized almost contact manifolds  $M_1 \times M_2$ could be generated by families of generalized almost contact  structures on the factor manifolds satisfying the conditions of Theorem \ref{T1}  if they were integrable with respect to other suitable derived brackets.  This study is left to a forthcoming paper.

Our final remark is on how the definitions of generalized co-K\"ahler and generalized Sasakian might be reformulated in terms of commuting pairs of generalized almost contact metric structures given Theorem \ref{T1}, which provides a more general organizing principle.  A definition of generalized co-K\"ahler structures was proposed by the author and Gomez in \cite{GT2}.  Specifically, we defined a generalized co-K\"ahler structure to be a normal generalized contact metric structure $(M, \Phi, E_{\pm}, G)$ where $G\Phi$
is also strong with respect to the  Courant bracket.  A definition equivalent to this can be stated in terms of commuting pairs of generalized almost contact structures:
\begin{definition}
A generalized co-K\"ahler structure  on a manifold $M$ is a commuting pair of almost generalized contact metric structures $( M, \Phi, E_{\pm}, G)$ and $(M, \tilde{\Phi}, \tilde{E}_{\pm},$ $ G)$ such that $\tilde {\Phi} = G\Phi$, $\tilde{E}_+ = GE_+= E_-$, $\tilde{E}_- = GE_-= E_+$ and $\Phi$ and $\tilde{\Phi} = G\Phi$ are integrable with respect to the Courant bracket.
\end{definition}
We remark that this definition is consistent with the notion of a binormal (2,1) generalized almost contact structure in (\cite{V3}). 

In light of our analysis above of the classical Sasakian case, we propose that a definition of
a generalized Sasakian structure on a manifold $M$ should be in terms of families of commuting pair of generalized almost contact metric structures $( M, \Phi, E_{\pm}, G)$ and $(M, \tilde{\Phi}, \tilde{E}_{\pm}, G)$ satisfying the condition $E_{\pm }= \tilde{E}_{\pm}$ or the condition $E_{\pm } = \tilde{E}_{\mp}$ and such that $\Phi$, and  $\tilde{\Phi}$ are integrable with respect to a derived bracket.  There have been previous definitions proposed for the notion of generalized Sasakian by Vaisman \cite{V1,V2, V3}, Sekiya \cite{S}, Inglesias-Ponte and Wade \cite{IW}, and Wright \cite{W} all specific to the situation  $M\times \mathbb{R}^+$.  Vasiman and Sekiya define Sasakian in terms of integrability of generalized almost complex structures on $M \times R^+$ as opposed to defining the concept in terms of the intrinsic data on $M$.   The remark that that $\mathcal{J}_2$ is also equal to $RG\mathcal{J}_1R^{-1}$ in Example \ref{ex18} shows that these definitions would be consistent with a definition of generalized Sasakian defined in terms of commuting pairs closed under derived brackets.  Wright defines a generalized Sasakian structure as commuting pairs of generalized contact structures on $M$ with the same conditions on the sections $E_\pm$ by starting with generalized K\"ahler structure on spaces $M\times \mathbb{R}^+$ considering reductions of these structures to M.  This approach is strong in highlighting  the centrality of commuting pairs of generalized almost contact structures, but misses the subtlety of the integrability conditions and does not distinguish between generalized contact structures arising from classical co-K\"ahler  structures and  those arising from classical Sasakian structures.  Inglesias-Ponte and Wade defined generalized Sasakian as commuting pairs of generalized almost contact structures ( which for them are $\mathcal{E}^1(M)$ Dirac structures) that are closed under the Wade bracket. This may not be the full picture and is relevant only to generalized K\"ahler structures on $ M \times \mathbb{R}^+$.  For clarity and flexibility of computation it's clear that definitions of generalized coK\"aher and generalized Sasakian should be in terms of commuting pairs of generalized almost contact metric structures.  The above discussion shows that a notion of generalized Sasakian can never be formulated in terms of products of two Courant integrable generalized contact structures. The necessary step seems to be a complete understanding of when the integrability of families of commuting pairs of generalized almost contact metric structures with respect to a derived bracket implies the corresponding commuting pairs of general almost complex structures is integrable under the Courant bracket.

\section*{Acknowledgements}
The author acknowledges, with gratitude,  support from a Eugene M. Lang Faculty Fellowship from Swarthmore College and the  hospitality of the Department of Mathematics of Barnard College/ Columbia University during the time when this paper was written. The author would like to thank Ralph Gomez for ongoing discussions of this topic and helpful comments on drafts of this paper.


Janet Talvaccchia \\
Department of Mathematics \& Statistics \\
Swarthmore College\\
Swarthmore, PA 19081 \\
USA\\
{\it E-mail address}: {\tt jtalvac1@swarthmore.edu}\\[0.3cm]

\label{last}
\end{document}